\documentclass[12pt]{amsart}
\usepackage{amsmath, amssymb, epsfig}

\newlength{\algorithmwidth}
\theoremstyle{plain}
\newtheorem{theorem}{Theorem}[section]

\newtheorem{lemma}[theorem]{Lemma}

\theoremstyle{definition}

\theoremstyle{remark}
\newtheorem*{remark}{Remark}

\numberwithin{equation}{section}

\def \R {\mathbb{R}}

\def \E {\mathbb{E}}

\def \< {\langle}
\def \> {\rangle}
\def \^ {\widehat}

\newcommand{\pr}[2]{\langle {#1} , {#2} \rangle}

\newcommand{\defby}{\overset{\mathrm{\scriptscriptstyle{def}}}{=}}
\newcommand{\bigO}{\mathrm{O}}

\begin{document}
\bibliographystyle{plain}
\title{Randomized Kaczmarz solver for noisy linear systems}

\author{Deanna Needell}
\thanks{D.N.~is with the Dept.~of Statistics, Stanford University, 390 Serra Mall, Stanford CA 94305, USA. e-mail:
\texttt{dneedell@stanford.edu}.} %


\date{Jan. 18, 2009}

\begin{abstract}
  The Kaczmarz method is an iterative algorithm for solving systems of linear equations $Ax=b$.  Theoretical convergence rates for this algorithm were largely unknown until recently when work was done on a randomized version of the algorithm.  It was proved that for overdetermined systems, the randomized Kaczmarz method converges with expected exponential rate, independent of the number of equations in the system.  Here we analyze the case where the system $Ax=b$ is corrupted by noise, so we consider the system $Ax \approx b + r$ where $r$ is an arbitrary error vector. We prove that in this noisy version, the randomized method reaches an error threshold dependent on the matrix $A$ with the same rate as in the error-free case. We provide examples showing our results are sharp in the general context.
\end{abstract}
\maketitle
\section{Introduction}

The Kaczmarz method~\cite{K37:Angena} is one of the most popular solvers of overdetermined linear systems and has numerous applications from computer tomography to image processing.  It is an iterative method, and so therefore is practical in the realm of very large systems of equations.  The algorithm consists of a series of alternating projections, and is often considered a type of \textit{Projection on Convex Sets} (POCS) method.  Given a consistent system of linear equations of the form 
$$
Ax = b,
$$
the Kaczmarz method iteratively projects onto the solution spaces of each equation in the system.  That is, if $a_1, \ldots, a_m \in \R^n$ denote the rows of $A$, the method cyclically projects the current estimate orthogonally onto the hyperplanes consisting of solutions to $\pr{a_i}{x} = b_i$.  Each iteration consists of a single orthogonal projection.  The algorithm can thus be described using the recursive relation,
$$
x_{k+1} = x_k + \frac{b_i - \pr{a_i}{x_k}}{\|a_i\|_2^2}a_i,
$$
where $x_k$ is the $k^{th}$ iterate and $i = (k$ mod $m) + 1$.  

Although the Kaczmarz method is popular in practice, theoretical results on the convergence rate of the method have been difficult to obtain.  Most known estimates depend on properties of the matrix $A$ which may be time consuming to compute, and are not easily comparable to those of other iterative methods (see e.g. ~\cite{DH97:Therate},~\cite{G05:Onthe},~\cite{HN90:Onthe}).  Since the Kaczmarz method cycles through the rows of $A$ sequentially, its convergence rate depends on the order of the rows.  Intuition tells us that the order of the rows of $A$ does not change the difficulty level of the system as a whole, so one would hope for results independent of the ordering.  One natural way to overcome this is to use the rows of $A$ in a random order, rather than sequentially.  Several observations were made on the improvements of this randomized version~\cite{N86:TheMath,HM93:Algebraic}, but only recently have theoretical results been obtained~\cite{SV06:Arandom,SV09:Arand}.  


\subsection{Randomized Kaczmarz}
In designing a random version of the Kaczmarz method, it is necessary to set the probability of each row being selected. Strohmer and Vershynin propose in~\cite{SV06:Arandom,SV09:Arand} to set the probability proportional to the Euclidean norm of the row. Their revised algorithm can then be described by the following:
$$
x_{k+1} = x_k + \frac{b_{p(i)} - \pr{a_{p(i)}}{x_k}}{\|a_{p(i)}\|_2^2}a_{p(i)},
$$
where $p(i)$ takes values in $\{1, \ldots, m\}$ with probabilities $\frac{\|a_{p(i)}\|_2^2}{\|A\|_F^2}$.  Here and throughout, $\|A\|_F$ denotes the Frobenius norm of $A$ and $\|\cdot\|_2$ denotes the usual Euclidean norm or spectral norm for vectors or matrices, respectively.  We note here that of course, one needs some knowledge of the norm of the rows of $A$ in this version of the algorithm. In general, this computation takes $\bigO(mn)$ time.  However, in many cases such as the case in which $A$ contains Gaussian entries, this may be approximately or exactly known.  

In~\cite{SV06:Arandom,SV09:Arand}, Strohmer and Vershynin prove the following exponential bound on the expected rate of convergence for the randomized Kaczmarz method,
\begin{equation}\label{SVbound}
\mathbb{E}\|x_k - x\|_2^2 \leq \Big(1 - \frac{1}{R}\Big)^k\|x_0 - x\|_2^2,
\end{equation}
where $R = \|A^{-1}\|^2\|A\|_F^2$, $x_0$ is an arbitrary initial estimate, and $\E$ denotes the expectation (over the choice of the rows). Here and throughout, we will assume that $A$ has full column rank so that $\|A^{-1}\| \defby \inf\{M : M\|Ax\|_2 \geq \|x\|_2$ for all $x\}$ is well defined.  We comment here that this particular mixed condition number comes as an immediate consequence of the simple probabilities used within the randomized algorithm.  

The first remarkable note about this result is that it is essentially independent of the number $m$ of equations in the system.  Indeed, by the definition of $R$, $R$ is proportional to $n$ within a square factor of $\kappa(A)$, the condition number of $A$ ($\kappa(A)$ is defined as the ratio of the largest to smallest singular values of $A$).  This bound also demonstrates, however, that the Kaczmarz method is an efficient alternative to other methods only when the condition number is very small.  If this is not the case, then other alternative methods may offer improvements in practice.  

The bound~\eqref{SVbound} and the relationship of $R$ to $n$ shows that the estimate $x_k$ converges exponentially fast to the solution in just $\bigO(n)$ iterations.  Since each iteration requires $\bigO(n)$ time, the method overall has a $\bigO(n^2)$ runtime.  Being an iterative algorithm, it is clear that the randomized Kaczmarz method is competitive only for very large systems.  For such large systems, the runtime of $\bigO(n^2)$ is clearly superior to, for example, Gaussian elimination which has a runtime of $\bigO(mn^2)$.  Also, since the algorithm needs only access to the randomly chosen rows of $A$, the method need not know the entire matrix $A$, which for very large systems is a clear advantage.  Thus the interesting cases for the randomized method are those in which $n$ and $m$ are large, and especially those in which $m$ is extremely large.  Strohmer and Vershynin discuss in detail in Section 4.2 of~\cite{SV09:Arand} cases where the randomized Kaczmarz method even outperforms the conjugate gradient method (CGLS). They show that for example, randomized Kaczmarz computationally outperforms CGLS for Gaussian matrices when $m > 3n$.  Numerical experiments in~\cite{SV09:Arand} also demonstrate advantages of the randomized Kaczmarz method in many cases.


Since the results of \cite{SV06:Arandom,SV09:Arand}, there has been some further discussion about the benefits of this randomized version of the Kaczmarz method (see \cite{CHJ09:Anote,SV09:Comments}).  The Kaczmarz method has been studied for over seventy years, and is useful in many applications.  The notion of selecting the rows randomly in the method has been proposed before (see \cite{N86:TheMath,CFMSS92,HM93:Algebraic}), and improvements over the standard method were observed. However, the work by Strohmer and Vershynin in~\cite{SV06:Arandom,SV09:Arand} provides the first proof on the rate of convergence. The rate is exponential in expectation and is in terms of standard matrix properties.  We are not aware of any other Kaczmarz method that provably achieves exponential convergence.

It is important to note that the method of row selection proposed in this version of the randomized Kaczmarz method is \textit{not} optimal, and an example that demonstrates this is given in~\cite{SV09:Arand}.  However, under this selection strategy, the convergence rates proven in~\cite{SV06:Arandom,SV09:Arand} are optimal, and there are matrices that satisfy the proven bounds exactly.  The selection strategy in this method was chosen because it often yields very good results, allows a provable guarantee of exponential convergence, and is computationally efficient.  

Since the algorithm selects rows based on their row norms, it is natural to ask whether one can simply scale the rows any way one wishes.  Indeed, choosing the rows based on their norms is related to the notion of applying a diagonal preconditioner.  However, since finding the optimal diagonal preconditioner for a system $Ax=b$ is itself a task that is often more costly than inverting the entire matrix, we select an easier, although not optimal, preconditioner that simply scales by the (square of the) row norms.  This type of preconditioner yields a balance of computational cost and optimality (see \cite{vdS69:Cond,S69:Optimality}).  The distinction between the effect of an alternative diagonal preconditioner on the Kaczmarz method versus the randomized method discussed here is important.  If the system is multiplied by a diagonal matrix, the standard Kaczmarz method will not change, since the angles between all rows do not change.  However, such a multiplication to the system in our randomized setting changes the probabilities of selecting the rows (by definition).  It is then not a surprise that this will also affect the convergence rate proved for this method (since multiplication will affect the value of $R$ in~\eqref{SVbound}).  

This randomized version of the Kaczmarz method provides clear advantages over the standard method in many cases.  Using the selection strategy above, Strohmer and Vershynin were able to provide a proof for the expected rate of convergence that shows exponential convergence.  No such convergence rate for any Kaczmarz method has been proven before.  These benefits lead one to question whether the method works in the more realistic case where the system is corrupted by noise.  In this paper we provide theoretical and empirical results to suggest that in this noisy case the method converges exponentially to the solution within a specified error bound.  The error bound is proportional to $\sqrt{R}$, and we also provide a simple example showing this bound is sharp in the general setting.

\section{Main Results}

Theoretical and empirical studies have shown the randomized Kaczmarz algorithm to provide very promising results.  Here we show that it also performs well in the case where the system is corrupted with noise.  In this section we consider the consistent system $Ax=b$ after an error vector $r$ is added to the right side:
$$
Ax \approx b + r.
$$
Note that we do not require the perturbed system to be consistent.  First we present a simple example to gain intuition about how drastically the noise can affect the system.  To that end, let $A$ be the $n \times n$ identity matrix, $b=0$, and suppose the error is the vector whose entries are all one, $r = (1, 1, \ldots, 1)$.  Then the solution to the noisy system is clearly $x = r = (1, 1, \ldots, 1)$, and the solution to the unperturbed problem is $x=0$.  
By Jensen's inequality, we have  
$$
\Big(\mathbb{E}\|x_k - r\|_2\Big)^2 \leq \mathbb{E}\Big(\|x_k - r\|_2^2\Big).
$$
Now considering the noisy problem, we may substitute $r$ for $x$ in~\eqref{SVbound}. Combining this with Jensen's inequality above, we obtain
\begin{equation}\label{X1}
\mathbb{E}\|x_k - r\|_2 \leq \Big(1 - \frac{1}{R}\Big)^{k/2}\|x_0-r\|_2.
\end{equation}
Then by the triangle inequality, we have
$$
\|r-x\|_2 \leq \|r-x_k\|_2 + \|x_k-x\|_2.
$$
Next, by taking expectation and using~\eqref{X1} above, we have
$$
\mathbb{E}\|x_k - x\|_2 \geq \|r - x\|_2 - \Big(1 - \frac{1}{R}\Big)^{k/2}\|x_0-r\|_2.
$$
Finally by the definition of $r$ and $R$, this implies
$$
\mathbb{E}\|x_k - x\|_2 \geq \sqrt{R} - \Big(1 - \frac{1}{R}\Big)^{k/2}\|x_0-r\|_2.
$$
This means that the limiting error between the iterates $x_k$ and the original solution $x$ is $\sqrt{R}$.  In~\cite{SV06:Arandom,SV09:Arand} it is shown that the bound provided in~\eqref{SVbound} is optimal, so even this trivial example demonstrates that if we wish to maintain a general setting, the best error bound for the noisy case we can hope for is proportional to $\sqrt{R}$.  Our main result proves this exact theoretical bound.

\begin{theorem}[Noisy randomized Kaczmarz]\label{thm}
Let $A$ have full column rank and assume the system $Ax=b$ is consistent. Let $x_k^*$ be the $k^{th}$ iterate of the noisy randomized Kaczmarz method run with $Ax \approx b + r$, and let $a_1, \ldots a_m$ denote the rows of $A$. Then we have
$$
\mathbb{E}\|x_k^* - x\|_2 \leq \Big(1 - \frac{1}{R}\Big)^{k/2}\|x_{0}\|_2 + \sqrt{R}\gamma,
$$
where $R = \|A^{-1}\|^2\|A\|_F^2$, $\gamma = \max_i \frac{|r_i|}{\|a_i\|_2}$, and the expectation is taken over the choice of the rows in the algorithm.
\end{theorem}

\begin{remark}
In the case discussed above, note that we have $\gamma = 1$, so the example indeed shows the bound is sharp.
\end{remark}

One may also recall the bound from perturbation theory (see e.g.~\cite{HJ85:Matrix-Analysis}) on the relative error in the perturbed case.  If we let $\hat{x} = A^\dagger (x+r)$ (where $A^\dagger \defby (A^*A)^{-1}A^*$ denotes the left inverse of $A$), then
$$
\frac{\|x-\hat{x}\|_2}{\|x\|_2} \leq \kappa(A) \frac{\|r\|_2}{\|Ax\|_2}.
$$
By applying the bound $\sqrt{R} \leq \kappa(A)\sqrt{n}$ to Theorem~\ref{thm} above, we obtain the bound
$$
\frac{\|x-\hat{x}\|_2}{\|x\|_2} \leq \kappa(A) \max_i\frac{\sqrt{n}|r_i|}{\|a_i\|_2\|x\|_2}.
$$
These bounds look similar in spirit, providing some more reassurance to the sharpness of the error bound.  It is important to note though that the first is obtained by applying the left inverse rather than an iterative method, which explains why the bounds are not exactly equal.  Of course for problems of large sizes, applying the inverse may not even be computationally feasible.

Before proving the theorem, it is important to first analyze what happens to the solution spaces of the original equations $Ax=b$ when the error vector is added.  Letting $a_1, \ldots a_m$ denote the rows of $A$, we have that each solution space $\pr{a_i}{x} = b_i$ of the original system is a hyperplane whose normal is $\frac{a_i}{\|a_i\|_2}$.  When noise is added, each hyperplane is translated in the direction of $a_i$.  Thus the new geometry consists of hyperplanes parallel to those in the noiseless case.  A simple computation provides the following lemma which specifies exactly how far each hyperplane is shifted.

\begin{lemma}\label{easylem}
Let $H_i$ be the affine subspaces of $\R^n$ consisting of the solutions to the unperturbed equations, $H_i = \{x: \left\langle a_i, x \right\rangle = b_i\}$. Let $H_i^*$ be the solution spaces of the noisy equations, $H_i^* = \{x: \left\langle a_i, x \right\rangle = b_i + r_i\}$. Then $H_i^* = \{w + \alpha_i a_i : w\in H_i\}$ where $\alpha_i = \frac{r_i}{\|a_i\|_2^2}$.
\end{lemma}
\begin{remark}
Note that this lemma does not imply that the noisy system is consistent. By definition of $H_i^*$ it is clear that each subspace is non-empty, but we are not requiring that the intersection of all $H_i^*$ be non-empty.
\end{remark}
\begin{proof}
First, if $w\in H_i$ then $\pr{a_i}{w + \alpha a_i} = \pr{a_i}{w} + \alpha\|a_i\|_2^2 = b_i + r_i$, so $w + \alpha a_i \in H_i^*$. Next let $u \in H_i^*$. Set $ w = u - \alpha a_i$. Then $\pr{a_i}{w} = \pr{a_i}{u} - r_i = b_i + r_i - r_i = b_i$, so $w \in H_i^*$. This completes the proof.
\end{proof}

We will also utilize the following lemma which is proved in the proof of Theorem 2 in~\cite{SV06:Arandom,SV09:Arand}.

\begin{lemma}\label{SVlem}
Let $x_{k-1}^*$ be any vector in $\mathbb{R}^n$ and let $x_k$ be its orthogonal projection onto a random solution space as in the noiseless randomized Kaczmarz method run with $Ax=b$. Then we have
$$
\mathbb{E}\|x_k - x\|_2^2 \leq \Big(1 - \frac{1}{R}\Big)\|x_{k-1}^* - x\|_2^2,
$$
where $R = \|A^{-1}\|^2\|A\|_F^2$, and the expectation is taken over the choice of the rows in the algorithm.
\end{lemma}

We are now prepared to prove Theorem~\ref{thm}.

\begin{proof}[of Theorem~\ref{thm}]
Let $x_{k-1}^*$ denote the $(k-1)^{th}$ iterate of noisy randomized Kaczmarz.  Using notation as in Lemma~\ref{easylem}, let $H_i^*$ be the solution space chosen in the $k^{th}$ iteration.  Then $x_k^*$ is the orthogonal projection of $x_{k-1}^*$ onto $H_i^*$. Let $x_k$ denote the orthogonal projection of $x_{k-1}^*$ onto $H_i$ (see Figure~\ref{fig0}). 

\begin{figure}[ht] 
  \includegraphics[scale=0.5]{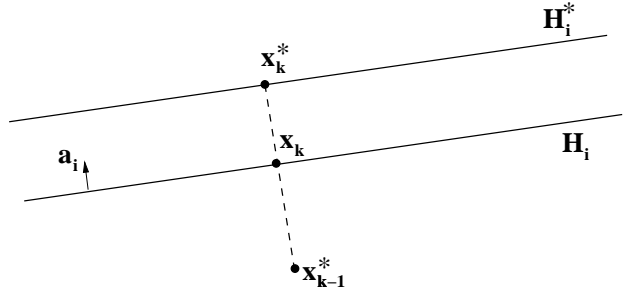}
  \caption{The parallel hyperplanes $H_i$ and $H_i^*$ along with the two projected vectors $x_k$ and $x_k^*$.}\label{fig0}
\end{figure}

By Lemma~\ref{easylem}  and the fact that $a_i$ is orthogonal to $H_i$ and $H_i^*$, we have that $x_k^* - x = x_k - x + \alpha_i a_i$.  Again by orthogonality, we have $\|x_k^* - x\|_2^2 = \|x_k - x\|_2^2 + \|\alpha_i a_i\|_2^2$.  Then by Lemma~\ref{SVlem} and the definition of $\gamma$, we have 
$$
\mathbb{E}\|x_k^* - x\|_2^2 \leq \Big(1 - \frac{1}{R}\Big)\|x_{k-1}^* - x\|_2^2 + \gamma^2,
$$
where the expectation is conditioned upon the choice of the random selections in the first $k-1$ iterations.
Then applying this recursive relation iteratively and taking full expectation, we have
\begin{align*}
\mathbb{E}\|x_k^* - x\|_2^2 &\leq \Big(1 - \frac{1}{R}\Big)^k\|x_{0} - x\|_2^2 + \sum_{j=0}^{k-1}\Big(1 - \frac{1}{R}\Big)^j\gamma^2\\
&\leq \Big(1 - \frac{1}{R}\Big)^k\|x_{0} - x\|_2^2 + R\gamma^2.
\end{align*}
By Jensen's inequality we then have
$$
\mathbb{E}\|x_k^* - x\|_2 \leq \left(\Big(1 - \frac{1}{R}\Big)^k\|x_{0} - x\|_2^2 + R\gamma^2\right)^{1/2} \leq \Big(1 - \frac{1}{R}\Big)^{k/2}\|x_{0} - x\|_2 + \sqrt{R}\gamma.
$$
This completes the proof.
\end{proof}

\section{Numerical Examples} 
 
\begin{figure}[h!]
\begin{center}
$\begin{array}{c@{\hspace{.1in}}c}
\includegraphics[width=2.5in]{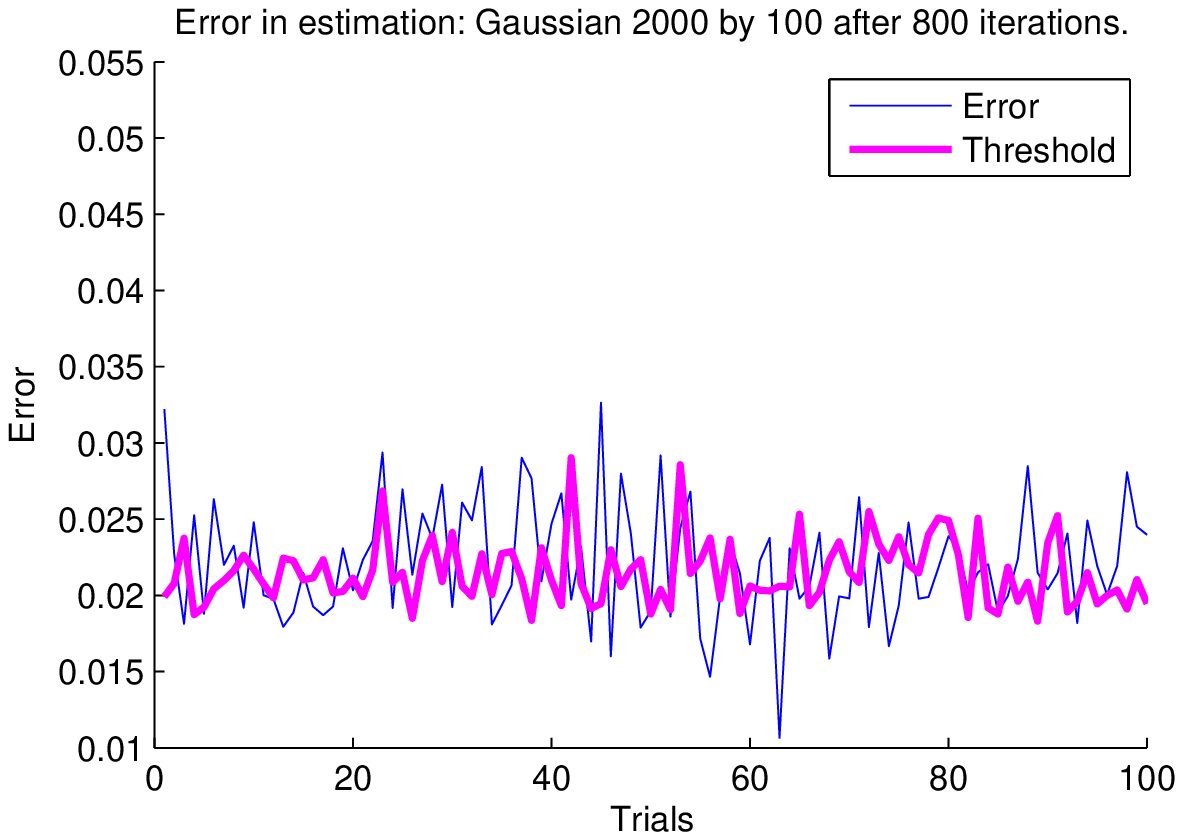}  &  
\includegraphics[width=2.5in]{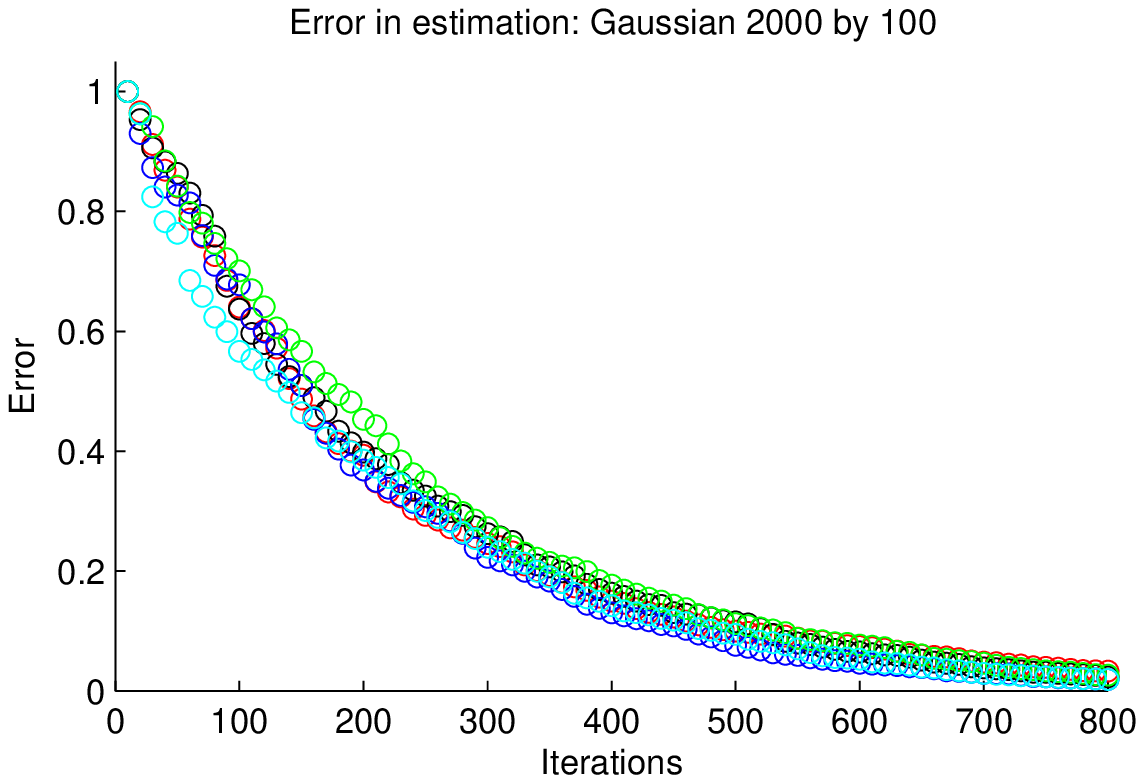}  \\ 
\includegraphics[width=2.5in]{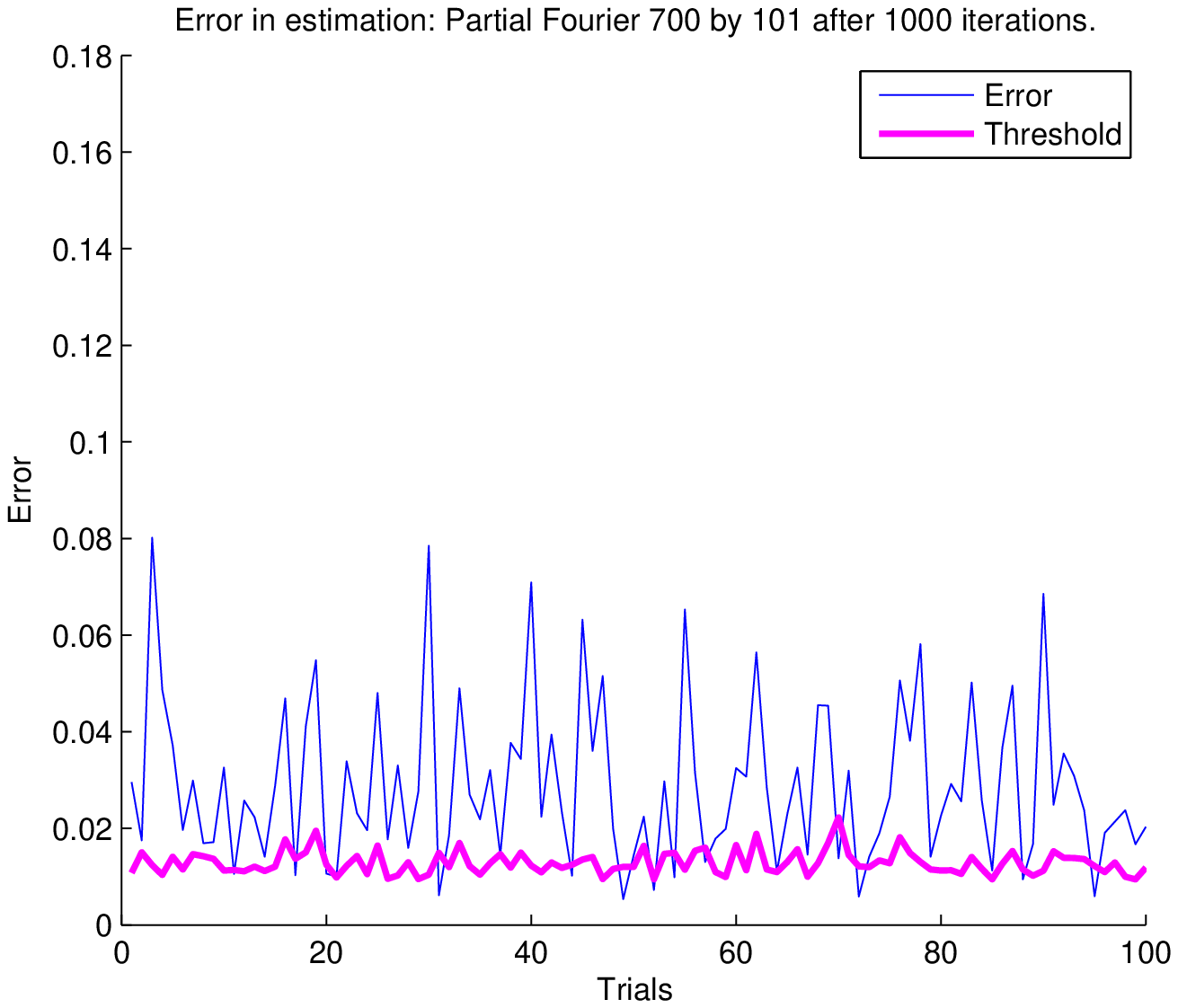}   &
\includegraphics[width=2.5in]{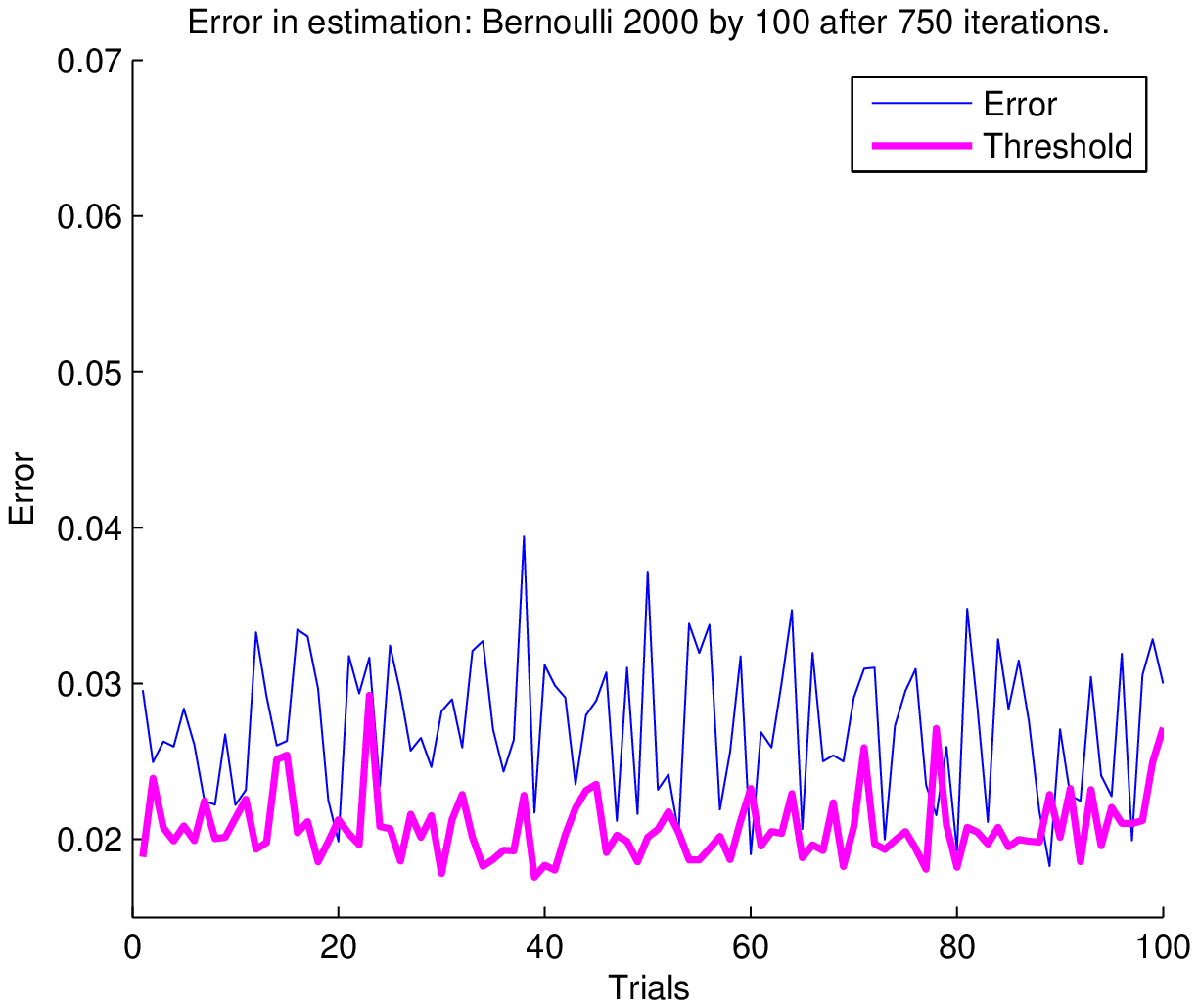} \\
\end{array}$
\end{center}
\caption{The comparison between the actual error in the randomized Kaczmarz estimate (thin line) and the predicted threshold (thick line). The mean values of $R$ in these experiments were 163.2 (upper left), 428.6 (lower left) and 162.4 (lower right). The scatter plot shows exponential convergence over several trials.}
\label{trio}
\end{figure}

In this section we describe some of our numerical results for the randomized Kaczmarz method in the case of noisy systems.  Figure~\ref{trio} depicts the error between the estimate by randomized Kaczmarz and the actual signal, in comparison with the predicted threshold value for several types of matrices.  The first study was conducted for 100 trials using $2000 \times 100$ Gaussian matrices (matrices who entries are i.i.d. Gaussian with mean $0$ and variance $1$) and independent Gaussian noise of norm $0.02$.  The systems were homogeneous, meaning $x=0$ and $b=0$.  The thick line is a plot of the threshold value, $\gamma\sqrt{R}$ for each trial.  The thin line is a plot of the error in the estimate after the given amount of iterations for the corresponding trial.  The scatter plot displays the convergence of the method over several randomly chosen trials from this study, and clearly shows exponential convergence.  The second study is a similar study but for the experiments in which we used partial Fourier matrices.  In this case we use $m=700$ and $n=101$.  For $j=1\ldots 700$ and $k=-50\ldots 50$, we set $A_{j,k} = \exp(2\pi ikt_j)$, where $t_j$ are generated uniformly at random on $[0, 1]$.  This type of generation is used to create nonuniformly spaced sampling values, and is used in many applications in signal processing, such as in the reconstruction of bandlimited signals.  The third study is similar but used matrices whose entries are Bernoulli ($0/1$ each with probability $0.5$).  All of these experiments were conducted to demonstrate that the error found in practice is close to that predicted by the theoretical results.  As is evident by the plots, the error is quite close to the threshold in all cases. 


%
%

 
 \subsection*{Acknowledgment}
I would like to thank Roman Vershynin for suggestions that simplified the proofs and for many thoughtful discussions. I would also like to thank Thomas Strohmer for his very appreciated guidance.
 
 \newpage
\bibliography{rk}

\end{document}